\documentclass[12pt]{article}
\usepackage[a4paper,margin=3cm]{geometry}
\usepackage{amsmath}
\usepackage{amssymb}

\newtheorem{theorem}{Theorem}[section]
\newtheorem{lemma}[theorem]{Lemma}

\newtheorem{proposition}[theorem]{Proposition}

\newenvironment{proof}{{\bf Proof.}}{\hfill$\Box$\\}
\newenvironment{definition}{{\vskip 3ex\bf Definition.} }{\\}
\newenvironment{remark}{{\vskip 3ex\bf Remark.} }{\\}

\newcommand{\N}{\mathbb{N}}

\newcommand{\R}{\mathbb{R}}

\newcommand{\ad}{\mathrm{ad}}
\newcommand{\rk}{\mathrm{rk}}
\newcommand{\tr}{\mathrm{tr}}
\newcommand{\Id}{\mathrm{Id}}

\newcommand{\spn}{\mathrm{span}}

\newcommand{\VV}{\mathcal{V}}
\newcommand{\WW}{\mathcal{W}}
\newcommand{\XX}{\mathcal{X}}
\newcommand{\YY}{\mathcal{Y}}

\newcommand{\CC}{\mathcal{C}}
\newcommand{\BG}{{\bf G}}
\newcommand{\BF}{{\bf F}}
\newcommand{\BV}{{\bf V}}
\newcommand{\BW}{{\bf W}}
\newcommand{\BX}{{\bf X}}
\newcommand{\BY}{{\bf Y}}
\newcommand{\BH}{{\bf H}}
\newcommand{\vp}{\varphi}

\newcommand{\fl}{\mathfrak l}

\newcommand{\pgl}{\mathfrak{pgl}}

\title{On contact equivalence of systems of ordinary differential equations}
\author{Wojciech Kry\'nski\thanks{
Institute of Mathematics, Polish Academy of Sciences, ul.~\'Sniadeckich 8, 00-956 Warszawa, Poland, E-mail:
krynski@impan.gov.pl}}

\begin{document}
\maketitle

\begin{abstract}
We consider a problem of equivalence of generic pairs $(\XX,\VV)$ on a manifold $M$, where $\VV$ is a distribution of rank $m$ and $\XX$ is a distribution of rank one. We construct a canonical bundle with a canonical frame. We prove that two pairs are equivalent if and only if the corresponding frames are diffeomorphic.

As a particular case, with $\VV$ integrable, we provide a new solution to the problem of contact equivalence of systems of $m$ ordinary differential equations: $x^{(k+1)}=F(t,x,x',\ldots,x^{(k)})$, where $k>2$ or $k=2$ and $m>1$. 
\end{abstract}



\section{Introduction}
Let $(F)$ be a system of $m$ ordinary differential equations given in the form
$$
x^{(k+1)}=F(t,x,x',\ldots,x^{(k)}),
$$
where $x\in\R^m$ and $F\colon\R^{1+(k+1)m}\to\R^m$. We consider equivalence problem up to the action of the group of contact transformations (or point transformations if $k=1$ and $m=1$). This problem was analysed by many authors. It was essentially solved by E.~J.~Wilczynski \cite{W}  in the linear case (see also Se-ashi \cite{YS} for modern approach). Recently, B.~Doubrov \cite{D,D2} extended Wilczynski invariants to the case of non-linear equations. He also characterized equations which are contact equivalent to the trivial one: $x^{(k+1)}=0$. Moreover, it is proved in the paper of B.~Doubrov, B.~Komarkov and T.~Morimoto \cite{DKM} that any system $(F)$ defines Cartan connection on a certain bundle over the space of $k$-jets and the problem of equivalence of equations is reduced to the problem of equivalence of Cartan connections.

Much is known about systems of low order. The problem was completely solved by E.~Cartan in the case of one equation of second order, i.e.~$k=1$, $m=1$. Systems of second order were analysed by M.~Fels \cite{F} and D.~Grossman \cite{G}. S-S.~Chern \cite{C} considered equations of order 3 and solved the equivalence problem via E.~Cartan method. He also found a class of equations which induce a conformal metric on the space of solutions. The class is characterized by the condition that all Wilczynski invariants vanish (so called Wuenschmann condition). A similar result was obtained by R.~Bryant \cite{B} in the order 4 case. Conformal structures associated to equations of arbitrary order were found by M.~Dunajski and P.~Tod \cite{DT} (see also \cite{K}).


The results of paper \cite{DKM} base on the observation that all contact information about $(F)$ is contained in the pair $(\XX,\VV)$, where $\XX$ and $\VV$ are two integrable distributions such that $\CC^k=\XX\oplus\VV$ is a canonical contact distribution on the space $J^k$ of $k$-jets. To be more precise $\VV$ is a vertical distribution tangent to the fibres of the projection $J^k\to J^{k-1}$ and $\XX$ is a line field spanned by the total derivative:
$$
X_F=\partial_t+\sum_{i=0}^{k-1}\sum_{j=1}^m x^{i+1}_j\partial_{x^i_j}+\sum_{j=1}^m F_j\partial_{x^k_j}.
$$
It is proved in \cite{DKM} that contact graded Lie algebra (i.e. symbol algebra of distribution $\CC^k$) together with the additional structure defined by the pair $(\XX,\VV)$ fulfils all assumptions of the general theorem of Morimoto \cite{M} and thus an arbitrary system $(F)$ defines, in a functorial way, a Cartan connection on a certain bundle.

In the present paper we give a different solution to the problem of contact equivalence of systems of ODEs. We assume that $k>2$ or $k=2$ and $m>1$. Our starting point is also the pair $(\XX, \VV)$. However we consider more general problem and in particular we drop the assumption that $\VV$ is integrable. This assumption was crucial in \cite{DKM} since it was important that the direct sum $\VV\oplus\XX$ has a fixed symbol algebra and due to this fact the general theory of graded Lie algebras was applied. In our approach the distribution $\VV$ can be a~priori arbitrary. Instead of integrability we impose additional conditions on the Lie brackets of sections of $\XX$ and $\VV$. Namely, we assume that the growth of dimensions of distributions $\ad^i_\XX\VV$ is maximal possible. In this way we define a class of \emph{regular} pairs $(\XX,\VV)$. This class is clearly generic in a sense that a small perturbation of a given pair $(\XX,\VV)$ is a regular pair. The set of all regular pairs contains as a subset the set of pairs which come from the differential equations (we say that $(\XX,\VV)$ is of \emph{equation type}). In next section we characterise pairs of \emph{equation type}. It appears that they satisfies two additional strong conditions.

The main result of the paper is the solution of the problem of equivalence of \emph{regular} pairs. For any $(\XX,\VV)$ which is \emph{regular} we construct the canonical bundle with the canonical frame and the problem of equivalence is reduced to the equivalence of such frames. As a by-product we obtain also a new solution to the problem of equivalence of systems of ODEs. It also appears that in this case our canonical frame in fact defines a Cartan connection. The main theorem is as follows:

\begin{theorem}\label{twierdzenie 0}
Let $M$ be a manifold of dimension $n=mk+m+1$, where $k>2$ or $k=2$ and $m>1$. For any regular pair $(\XX,\VV)$ on $M$ such that $\rk\,\VV=m$, there exists a canonical principal $\mathrm{PGL}(2)_0\oplus\mathrm{GL}(m)$-bundle $B(\XX,\VV)$ over $M$ which possesses a canonical frame. Two pairs $(\XX,\VV)$ and $(\YY,\WW)$ are equivalent if and only if the corresponding frames are diffeomorphic. The symmetry group of $(\XX,\VV)$ is at most $(n+m^2+2)$-dimensional and it has maximal dimension if and only if $(\XX,\VV)$ is locally of equation type and the corresponding system of equations is trivial: $x^{(k+1)}=0$.
\end{theorem}
The group $\mathrm{PGL}(2)_0$ is the group of real M\"obius transformations which preserve a fixed point in the projective line. It is consisted of matrices:
$$
\left(
\begin{array}{cc}
a&b\\
0&1\\
\end{array}
\right)\in\mathrm{GL}(2).
$$
If $t$ is a projective parameter on $\R P^1$, then the right action of $\mathrm{PGL}(2)_0$ on $t$ is defined in the following way:
$$
t\mapsto\frac{at}{bt+1},
$$
so that $0$ is mapped to $0$.

We have already mentioned that the frame in the main theorem can be viewed as a Cartan connection (if $(\XX,\VV)$ is of equation type). Indeed, we show in Section~5 that the dual coframe defines Cartan connection of type $(L,\mathrm{PGL}(2)_0\oplus\mathrm{GL}(m))$, where $L$ is a semi-direct product of $\mathrm{SL}(2)$ and the following affine-like Lie group:
$$
A=\left\{\left(
\begin{array}{cc}
\Id_{k+1}&0\\
V&G
\end{array}
\right)\in \mathrm{GL}(m+k+1)\ |\ V\in\R^{m\times(k+1)}, G\in \mathrm{GL}(m)\right\},
$$
where $\mathrm{SL}(2)$ acts irreducibly on each $\R^{k+1}$ being the component of $\R^{m\times(k+1)}$.

Let us stress that our main theorem holds for $k>2$ or $k=2$ and $m>1$. The dimension of maximal symmetry group in the cases $k=1$ or $k=2$ and $m=1$ do not fit the general scheme $(m^2+(k+1)m+3)$ (see \cite{DKM,F,O}). For example the group of symmetries of the equation $x'''=0$ is 10-dimensional, as it is proved in \cite{C}.

\section{Geometry of ODE}
Let us consider a system $(F)$ of $m$ ordinary differential equations of order $k+1$:
$$
x^{(k+1)}=F(t,x,x',\ldots,x^{(k)}),
$$
where $x\in\R^m$ and $F=(F_1,\ldots,F_m)$ is a function $\R^{1+(k+1)m}\to\R^m$. In this section we will focus on a geometric description of $(F)$ (we refer to \cite{Enc} for the details).

Recall that $J^{k+1}(1,m)$ denotes the space of $k+1$ jets of functions $\R\to\R^m$. Let $(t, x^0,\ldots, x^{k+1})$, where $x^i=(x^i_1,\ldots,x^i_m)$, be the standard coordinate system on $J^{k+1}(1,m)$. The system $(F)$ is equivalently defined by a certain corank $m$ sub-manifold $E_F\subseteq J^{k+1}(1,m)$. Namely:
$$
E_F=\{(t, x^0,\ldots, x^{k+1})\in J^{k+1}(1,m)\ |\ x^{k+1}_j-F_j(t,x^0,\ldots,x^k)=0,\ \ j=1,\ldots,m\}.
$$
Note that the functions $t, x^0,\ldots, x^k$ restricted to $E_F$ can be considered as coordinates, since $E_F$ projects regularly on the subspace $\{x^{k+1}=0\}\subseteq J^{k+1}(1,m)$. In particular $E_F$ and $J^k(1,m)$ are diffeomorphic in the canonical way. Let
$$
\omega^i_j=dx^i_j-x^{i+1}_jdt,
$$
$i=0,\ldots,k$, $j=1,\ldots,m$, denote \emph{contact 1-forms} on the space of jets.

Consider restrictions of $\omega^i_j$ to $E_F$. We define $\XX_F$ as the intersection of all $\ker\omega^i_j|_{E_F}$ for every $i,j$. Note that $\XX_F$ is rank-one distribution (line field) spanned by vector field $X_F$, called \emph{total derivative}. In the coordinates it takes the form:
$$
X_F=\partial_t+\sum_{i=0}^{k-1}\sum_{j=1}^m x^{i+1}_j\partial_{x^i_j}+\sum_{j=1}^m F_j\partial_{x^k_j}.
$$

\emph{Vertical distribution} $\VV_F$ on $E_F$ is defined as the kernel of the canonical projection $E_F\to J^{k-1}(1,m)$. In coordinates:
$$
\VV_F=\spn\{\partial_{x^k_j}\ |\ j=1,\ldots,m\}.
$$

In this way we assigned the pair:
$$
(\XX_F,\VV_F)
$$
to a given system $(F)$. Recall that \emph{contact distributions} are defined as:
$$
\CC^i=\bigcap_{s=0,\ldots,i}\ker\omega_1^s\cap\ldots\cap\ker\omega_m^s.
$$
We get the following decomposition:
$$
\CC^k=\XX_F\oplus\VV_F
$$

A \emph{contact transformation} is a mapping $\Psi\colon J^{k+1}(1,m)\to J^{k+1}(1,m)$ such that $\Psi_*(\CC^k)=\CC^k$. We say that two equations $(F)$ and $(G)$ are \emph{contact equivalent} if there exists a contact transformation such that $\Psi|_{E_F}$ is a diffeomorphism of $E_F$ onto $E_G$.

We have the following:
\begin{proposition}\label{stwierdzenie 0}
System $(F)$ and $(G)$ are equivalent if and only if there exists a diffeomorphism $\Phi\colon E_F\to E_G$ which transforms $(\XX_F,\VV_F)$ onto $(\XX_G,\VV_G)$.
\end{proposition}
\begin{proof}
See Theorem 1 \cite{DKM}. The Proposition is a consequence of Lie-B\"acklund theorem.
\end{proof}

Consider an arbitrary pair $(\XX,\VV)$ on a manifold $M$, where $\VV$ is rank $m$ distribution and $\XX$ is a line field. Assume $\dim M=\dim J^k(1,m)$. We will see that not every pair $(\XX,\VV)$ is equivalent to a pair defined by an equation. Therefore we introduce the following notion.
\begin{definition}
Let $\XX$ be a line field and $\VV$ be a distribution on a manifold $M$. The pair $(\XX,\VV)$ is of \emph{equation type} if there exists a system $(F)$ and a diffeomorphism $\Phi\colon M\to E_F$ such that
$\Phi_*(\XX)=\XX_F$ and $\Phi_*(\VV)=\VV_F$. The pair $(\XX,\VV)$ is \emph{locally of equation type} if for any $x\in M$ there exists a neighbourhood $U\ni x$ such that $(\XX|_U,\VV|_U)$ is of equation type.
\end{definition}

For a given distribution $\WW$ we denote by $\Gamma(\WW)$ the set of all smooth sections of $\WW$. If $X$ is a vector field then we denote by $\ad_X\WW$ the distribution spanned by all Lie brackets $[X,Y]$, where $Y\in\Gamma(\WW)$. Note that it may happen that $\ad_X\WW$ is not of constant rank, even if $\WW$ was. However we will assume that it is not the case. Locally, a distribution $\WW$ of constant rank $m$ can be written as:
$$
\WW=\spn\{Y_1,\ldots,Y_m\}
$$
for certain $Y_1,\ldots,Y_m\in\Gamma(\WW)$. Any such tuple $Y=(Y_1,\ldots,Y_m)$ will be called \emph{a local frame of $\WW$}.

Consider a pair $(\XX,\VV)$, where $\XX$ is a line field and $\VV$ is a distribution of rank $m$. Let $X\in\Gamma(\XX)$. For a pair $(X,\VV)$, we introduce a sequence of distributions $\VV^0_X\subseteq\VV^1_X\subseteq\ldots$ defined inductively as follows:
$$
\VV^0_X=\VV,\qquad \VV^{i+1}_X=\ad_X\VV^i_X.
$$
We also define the distributions:
$$
\VV^i=\VV_X^i+\XX.
$$
It is clear that $\VV^i$ are independent on the choice of a nowhere-vanishing section $X$ of $\XX$ (our considerations will be local, therefore we can assume that such a section $X$ exists).

From K.~Yamaguchi \cite{Y} we easily deduce the following characterisation of pairs $(\XX,\VV)$ which are of equation type.
\begin{theorem}\label{twierdzenie 1}
A pair $(\XX,\VV)$, where $\rk\,\VV=m$, on a manifold $M$ of dimension $n$ is locally of equation type if and only if there exists $k\in\N$ such that the following conditions are satisfied
\begin{itemize}
\item[(G1)]$\rk\VV^i=(i+1)m+1$, for $i=0,\ldots,k$,
\item[(G2)]$\rk\VV^k=n$,
\item[(G3)]there exists an integrable $\WW^i\subseteq\VV^i$ of corank 1, for every $i=0,\ldots,k$,
\item[(G4)]$\WW^0=\VV$ and $\WW^i=\mathrm{Ch}(\VV^{i+1})$, for $i=0,\ldots,k-1$.
\end{itemize}
\end{theorem}
In above the symbol $\mathrm{Ch}(\WW)$ stands for the Cauchy characteristic of $\WW$, i.e.
$$
\mathrm{Ch}(\WW)=\{Y\in\Gamma(\WW)\ |\ \ad_Y\WW\subseteq\WW\}.
$$

\section{Normal Form}
Our main aim is to solve the local equivalence problem for systems $(F)$. By Theorem~\ref{twierdzenie 1} we can consider pairs $(\XX,\VV)$ satisfying the conditions (G1)-(G4) instead of differential equations. However, in order to have a more general perspective we will not assume the integrability conditions (G3) and (G4).

\begin{definition}
A pair $(\XX,\VV)$ is \emph{regular} if conditions (G1) and (G2) hold.
\end{definition}

Note that a generic pair $(\XX,\VV)$ on $\R^{(k+1)m+1}$ is regular in a neighbourhood of $0$.

In the current section we will construct the canonical frames of $\XX$ and $\VV$. The procedure is an implementation of the concept of Laguerre-Forsyth normal form of~$(F)$, which was already used by Wilczynski \cite{W}.

In the sequel we will use the matrix notation. For a tuple of vector fields $V=(V_1,\ldots,V_m)$ we will denote $\ad_XV=(\ad_XV_1,\ldots,\ad_XV_m)$. If $G=(G^i_j)_{i,j=1,\ldots,m}$ is a $\mathrm{GL}(m)$-matrix valued function, then we will write $VG=(\sum_iV_iG_1^i,\ldots,\sum_iV_iG_m^i)$. Note that from (G1) and (G2) it follows that if $V$ is a local frame of $\VV$ and $X$ is a non-vanishing section of $\XX$ then $(X,V, \ad_XV,\ldots,\ad_X^kV)$ constitutes a local frame on the manifold $M$. The following lemma can be treated as a nonlinear version of the first step in the construction of the Laguerre-Forsyth normal form of~$(F)$.

\begin{lemma}\label{lemat 1}
Let $(\XX,\VV)$ be a regular pair. Then, for any non-vanishing section $X$ of $\XX$, there exists a local frame $V=(V_1,\ldots,V_m)$ of $\VV$, such that $\ad^{k+1}_XV_i=0\mod\VV^{k-1}$ for $i=1,\ldots,m$.
\end{lemma}
\begin{proof}
Let $W=(W_1,\ldots, W_m)$ be any local frame of $\VV$. We will find functions $G=(G^j_i)_{i,j=1,\ldots,m}$ such that $V=WG$ is the desired frame. In the matrix notation we have
$$
\ad^{k+1}_XV = (\ad^{k+1}_XW)G + (k+1)(\ad^k_XW)X(G) \mod \VV^{k-1}.
$$
Assume that
$$
\ad^{k+1}_XW = (\ad^k_XW)H \mod\VV^{k-1},
$$
for a certain $H=(H_i^j)_{i,j=1,\ldots,m}$. Since $\ad^k_XV=(\ad^k_XW)G\mod\VV^{k-1}$, we obtain the following equation
\begin{equation}\label{rownanie 1}
HG+(k+1)X(G)=0
\end{equation}
It can be solved locally and, if $G$ is a solution, then $\ad^{k+1}_XV=0\mod\VV^{k-1}$.
\end{proof}

\begin{definition}
A local frame $V=(V_1,\ldots,V_m)$ of $\VV$ is called a \emph{normal frame of $(X, \VV)$} if $\ad^{k+1}_XV=0\mod\VV^{k-1}$.
\end{definition}

Lemma \ref{lemat 1} implies that normal frames exist. They are solutions to the system \eqref{rownanie 1} of first ODEs. Let $V$ be a normal frame of a pair $(X,\VV)$. Then there exist matrix valued functions $K_0,\ldots, K_{k-1}$ defined by the equation:
$$
\ad_X^{k+1}V+(\ad_X^{k-1}V)K_{k-1}+\cdots+(\ad_XV)K_1+VK_0=0.
$$
Note that by definition of normal frame there is no term of $k$-th order in the equation above. If $W=VG$ is a different normal frame of $(X,\VV)$ then $X(G) = 0$ (see equation \eqref{rownanie 1} with $H\equiv 0$). Hence
$$
\ad^i_XW =(\ad^i_XV)G
$$
for every $i$ and we see that:
$$
\ad_X^{k+1}W+(\ad_X^{k-1}W)\tilde K_{k-1}+\cdots+(\ad_XW)\tilde K_1+W\tilde K_0=0
$$
where:
$$
\tilde K_i=G^{-1}K_iG.
$$
It follows that a pair $(X,\VV)$ defines operators:
$$
K_i\colon\VV\to\VV,
$$
which are vector bundle homorphisms.

\begin{remark}
Let $(F)$ be the geodesic equation
$$
x_i''=-\sum_{p,q=1}^m\Gamma^i_{pq}x_p'x_q',
$$
on $m$ dimensional manifold $M$, where $\Gamma=(\Gamma^i_{pq})$ is an affine connection. Without lost of generality we may assume that $\Gamma$ is symmetric, $\Gamma^i_{pq}=\Gamma^i_{qp}$. Consider the pair $(\XX_F,\VV_F)$. We can identify $J^1(1,m)\simeq TM\times\R$ and then $x^0=(x^0_1,\ldots,x^0_m)$ are coordinates on $M$ whereas $x^1=(x^1_1,\ldots,x^1_m)$ are the corresponding linear coordinates on the tangent spaces. Let
$$
X_F=\partial_t+\sum_sx^1_s\partial_{x^0_s} -\sum_{s,p,q}\Gamma^s_{pq}x_p^1x_q^1\partial_{x^1_s}
$$
be the total derivative. Then $K_0$, defined by the pair $(X_F,\VV_F)$, encodes the curvature tensor of $\Gamma$. Namely it can be proved that
$$
K_0=\left(\sum_{p,q}R^j_{ipq}x^1_px^1_q\right)_{i,j=1,\ldots,m},
$$
where $R^j_{ipq}=\partial_{x^0_i}(\Gamma^j_{pq}) -\partial_{x^0_q}(\Gamma^j_{ip})+\Gamma^j_{ir}\Gamma^r_{pq} -\Gamma^r_{ip}\Gamma^j_{rq}$.
\end{remark}

\begin{remark}
The general theory of pairs $(X,\VV)$ is developed in \cite{JK}. It appears that such pairs are very important in the theories of sprays, Veronese webs, affine control systems and ODEs with fixed time-scale.
\end{remark}

Unfortunately, the dependence of the operators $K_i$ on the choice of a section $X$ of $\XX$ is not tensorial. However, certain special sections can be chosen. Indeed, the next lemma can be regarded as the second step in the nonlinear version of the construction of the Laguerre-Forsyth normal form of $(F)$. We will write $K_i^X$ in order to distinguish operators corresponding to different sections $X$ of line bundle $\XX$.
\begin{lemma}\label{lemat 2}
$$
\tr\,K_{k-1}^{fX} = f^2\tr\,K_{k-1}^{X} - mc_kS^X(f),
$$
where $c_k=-\frac{1}{24}k(k+1)(k+2)$, and
$$
S^X(f)=2fX^2(f)-X(f)^2.
$$
\end{lemma}
\begin{proof}
Let $V$ be normal frame for $X$, and $VG$ be normal frame for $fX$. We compute directly:
\begin{eqnarray}\label{rownanie 2}
\ad^{k+1}_{fX}(VG)&=&(\ad^{k+1}_XV)f^{k+1}G+
\sum_{i=0}^k(\ad^k_XV)f^{k+1-i}X(f^iG)\\
&+&\sum_{i=0}^{k-1}\sum_{j=0}^{k-1-i}(\ad^{k-1}_XV)f^{k-i-j}X(f^{j+1}X(f^iG))
\mod\ad^{k-2}_X\VV,\XX.\nonumber
\end{eqnarray}
From normality of frames $V$ and $VG$ it follows $\sum_{i=0}^k f^{k+1-i}X(f^iG)=0.$ Applying the Leibnitz rule and dividing over $f^k$ simplifies this equation to:
\begin{equation}\label{rownanie 3}
fX(G)=-\frac{k}{2}X(f)G.
\end{equation}
Differentiating both sides with $fX$ and replacing $fX(G)$ with $-\frac{k}{2}X(f)G$ gives:
\begin{equation}\label{rownanie 4}
f^2X^2(G)=\left(\left(\frac{k}{2}+\frac{k^2}{4}\right)X(f)^2-
\frac{k}{2}fX^2(f)\right)G.
\end{equation}
Now, we can substitute \eqref{rownanie 3} and \eqref{rownanie 4} to equation \eqref{rownanie 2} and find the exact formula for the coefficient next to $\ad_X^{k-1}V$ in terms of $f$, $G$, $X(f)$ and $X^2(f)$. In this way the lemma follows (see \cite{K} for detailed computations).
\end{proof}
\begin{remark}
Operator $S^X$ is called \emph{Schwartzian}. It is justified by the following reasoning. Let $\gamma\colon t\mapsto\gamma(t)$ be a trajectory of $X$. Then $X(g)\circ\gamma = \frac{d}{dt}(g\circ\gamma)$ for any function $g\colon M\to\R$. Let us consider a reparametrization $\vp\colon s\mapsto\vp(s)$, such that $\gamma\circ\vp$ is a trajectory of the vector field $fX$. Then $f(\gamma\circ\vp)=\vp'$ and:
\begin{eqnarray*}
S^X(f)(\gamma\circ\vp)&=&2\vp'\left(\frac{1}{\vp'}\frac{d}{ds}\right)^2(\vp')- \left(\frac{1}{\vp'}\frac{d}{ds}\vp'\right)^2\\
&=&2\frac{d}{ds}\left(\frac{\vp''}{\vp'}\right)-\left(\frac{\vp''}{\vp'}\right)^2\\
&=&2\frac{\vp'''}{\vp'}-3\left(\frac{\vp''}{\vp'}\right)^2.
\end{eqnarray*}
The last term is equal to the Schwartz derivative of $\vp$.
\end{remark}

It is known that parameterizations corresponding to different solutions of M\"obius equation $S^X(f)=0$ are related by the formula:
$$
s=\frac{at+c}{bt+d},
$$
where $\left(\begin{array}{cc}a&b\\c&d\end{array}\right)\in GL(2)$.
Therefore, on each integral line of $\XX$ there exists the canonical projective structure. 

\begin{definition}
A section $X$ of $\XX$ is called a \emph{projective vector field of $(\XX,\VV)$} if $\tr\,K_{k-1}^X=0$.
\end{definition}

Lemma \ref{lemat 2} implies that the projective vector fields form a 2-parameter family on any integral curve of $\XX$. We get that the group
$$
\mathrm{PGL}(2)_0=\left\{\left(
\begin{array}{cc}
a&b\\
0&1\\
\end{array}
\right)\in GL(2)\right\}
$$
acts on the set of all projective vector fields along an integral line of $\XX$ freely and transitively. As a conclusion of the current section we summarise that if $X$ is a projective vector field and $V$ is a corresponding normal frame then
$$
\ad_X^{k+1}V+(\ad_X^{k-1}V)K_{k-1}^X+\cdots+(\ad_XV)K_1^X+VK_0^X=0,
$$
where $\tr\,K_{k-1}^X=0$.
\section{Canonical Frame}
We are now in the position to construct a canonical principal bundle for a regular pair $(\XX,\VV)$ on a manifold $M$. Let $x\in M$ and let $B(\XX,\VV)(x)$ be the set consisting of pairs $(\chi,\nu)$ where $\nu$ is a linear basis of $\VV(x)$ and $\chi$ is a germ at $x$ of a projective vector field along the integral line of $\XX$ which passes through $x$. In fact, since projective vector fields are described by the second order ODE, a germ $\chi$ depend only on the 1-jet at $x$ of a section $X\in\Gamma(\XX)$. In particular $B(\XX,\VV)(x)$ is a finite dimensional manifold. We will denote by $j^0\chi$ the 0-jet of $\chi$ at $x$. It follows that $j^0\chi$ is a vector in $\XX(x)$. We define
$$
B(\XX,\VV)=\bigcup_{x\in M}B(\XX,\VV)(x).
$$
Clearly $B(\XX,\VV)$ is principal $\mathrm{GL}(m)\oplus \mathrm{PGL}(2)_0$-bundle over $M$, as follows from the previous section. We call it \emph{the canonical bundle of $(\XX,\VV)$}. The projection $B(\XX,\VV)\to M$ will be denoted by $\pi$.

On $B(\XX,\VV)$ there are canonical fundamental vertical vector fields which come from the infinitesimal action of the structural group. Namely, vector fields $\BG^s_t$, where $s,t=1,\ldots,m$ come from the action of $\mathrm{GL}(m)$ and $\BF^0$ and $\BF^1$ come from the action of $\mathrm{PGL}(2)_0$. We will abbreviate $\BG=(\BG^s_t)_{s,t=1,\ldots,m}$ and $\BF=(\BF^0,\BF^1)$. If we choose a local horizontal section $M\ni x\mapsto(\chi(x),\nu(x))\in B(\XX,\VV)(x)$ then we can introduce local coordinates on fibres of $B(\XX,\VV)$. Indeed, any point in $B(\XX,\VV)$ can be represented as $(\chi(x),\nu(x))\cdot (F,G)$, where $x\in M$, $G=(G^s_t)\in\mathrm{GL}(m)$ and
$$
F=\left(\begin{array}{cc}
F_0&F_1\\
0&1\\
\end{array}\right)\in\mathrm{PGL}(2)_0.
$$
In this coordinates:
$$
\BG_q^p=\sum_{j=1}^mG^j_q\partial_{G^j_p},\qquad \BF^0=F_0\partial_{F_0},\qquad \BF^1=F_0\partial_{F_1}.
$$
Note that:
$$
\BG_q^p(G_p^r)=G_q^r,\qquad \BF^0(F_0)=F_0,\qquad \BF^1(F_1)=F_0,
$$
and there is no other nontrivial differentiation of functions $G,F_0,F_1$ in the directions of $\BG$ and $\BF$. The definition of $\BG$ and $\BF$ implies:
\begin{proposition}
The following structural equations are satisfied:
\begin{equation}\label{strukturalne 0}
[\BG^p_q,\BG^s_t]=\delta^p_t\BG^s_q-\delta^s_q\BG^p_t,\qquad[\BG,\BF]=0,\qquad[\BF^0,\BF^1]=\BF^1.
\end{equation}
\end{proposition}

Our aim is to choose additional vector fields: $\BX$ and $\BV^i_j$, $i=0,\ldots,k$, $j=1,\ldots,m$, such that the tuple 
$$
(\BG^s_t,\BF^0,\BF^1,\BX,\BV^0_j,\ldots,\BV^k_j\ |\ s,t,j=1,\ldots,m)
$$
constitutes a frame on $B(\XX,\VV)$. We will briefly write $\BV^i=(\BV^i_1,\ldots,\BV^i_m)$.

Vector field $\BX$ is defined by the following lemma.
\begin{lemma}\label{lemat 4}
There is the unique vector field $\BX$ on the canonical bundle $B(\XX,\VV)$ such that if $t\mapsto p(t)=(\chi(t),\nu(t))$ is an integral curve of $\BX$, then $X(t)=j^0\chi(t)=\pi_*(\BX(p(t))$ is a projective vector field on $M$ along $t\mapsto\pi(p(t))$, and $V(t)=\nu(t)$ is a normal frame of $(X,\VV)$.
\end{lemma}
\begin{proof}
Take an arbitrary projective vector field $X$ and some corresponding normal frame $V$ of $\VV$. Then, $(X,V)$ defines a horizontal section $\Gamma\subseteq B(\XX,\VV)$. We define $\BX$ on $\Gamma$ as the lift $\BX|_\Gamma=\pi^{-1}(X)$. The construction is correct, since, by \eqref{rownanie 1} and \eqref{rownanie 3}, $X$ and $V$ are uniquely determined along one integral curve of $\XX$ by the system of ODEs and the initial condition. The initial condition is given by a point in $B(\XX,\VV)$. In other words there is a partial connection on $B(\XX,\VV)$ in the direction of the line field $\XX$.
\end{proof} 

\begin{lemma}\label{lemat 5}
Let $X$ be a projective vector field and $V$ be a normal frame of $(X,\VV)$. In local coordinates on $B(\XX,\VV)$ defined by $X$ and $V$ we have:
$$
\BX=\frac{1}{F_0}X -2\frac{F_1}{F_0}\BF^0 -\frac{(F_1)^2}{(F_0)^2}\BF^1 -k\frac{F_1}{F_0}\sum_{j=1}^m\BG^j_j.
$$
\end{lemma}
\begin{proof}
Let us assume that $X=\partial_t$ is a projective vector field, and consider the action of $\mathrm{PGL}(2)_0$. Let
$$
s=\vp(t)=\frac{at}{bt+1}.
$$
Then $F_0(0)=a$ and $F_1(0)=b$ are coordinates on $B(\XX,\VV)$ at points where $t=s=0$ and our aim is to check how $F_0$ and $F_1$ change when $s=s_0\neq 0$. We have
$$
t=\vp^{-1}(s)=\frac{s}{a-bs}.
$$
Let $t_0=\vp^{-1}(s_0)$. We set $\tilde s=s-s_0$ and $\tilde t=t-t_0$. Then, by definition
$$
\tilde s=\frac{F_0(s_0)\tilde t}{F_1(s_0)\tilde t+1}.
$$
We compute
$$
\tilde s=s-s_0=\frac{at}{bt+1}-s_0= \frac{a(\tilde t+t_0)}{b(\tilde t+t_0)+1}-s_0= \frac{\left(\frac{a-bs_0}{bt_0+1}\right)\tilde t} {\left(\frac{b}{bt_0+1}\right)\tilde t+1},
$$
where we use the simple formula $at_0-s_0-bt_0s_0=0$. From above it follows that
$$
F_0(s_0)=\frac{a-bs_0}{bt_0+1},\qquad F_1(s_0)=\frac{b}{bt_0+1},
$$
and if we express $t_0$ in terms of $s_0$ we get
$$
F_0(s_0)=\frac{(a-bs_0)^2}{a},\qquad F_1(s_0)=\frac{b(a-bs_0)}{a}.
$$
Therefore
$$
\frac{d}{ds}F_0(s)=-2F_1(s),\qquad \frac{d}{ds}F_1(s)=-\frac{(F_1(s))^2}{F_0(s)},
$$
and additionally we have
$$
\frac{d}{ds}\vp^{-1}(s)=\frac{1}{F_0(s)}.
$$
Hence
$$
\BX=\frac{1}{F_0}X-2F_1\partial_{F_0} -\frac{(F_1)^2}{F_0}\partial_{F_1}\mod\partial_{G}.
$$
The coefficient next to $\partial_{G}$ remains unknown. However form equation \eqref{rownanie 3} it follows that the evolution of $G$ is given by the equation
$$
\frac{d}{ds}G^p_q(s)= -\frac{k}{2}G^p_q(s)F_0(s)\frac{d}{ds}\left(\frac{1}{F_0(s)}\right) =-k\frac{F_1(s)}{F_0(s)}G^p_q(s).
$$
We also have $\sum_{p,q}G^p_q\partial_{G^p_q}=\sum_j\BG^j_j$ and this formula completes the proof.
\end{proof}

The choice of $\BV^i_j$ is more complicated since there is no connection on $TM$ in the directions of $\VV$. Therefore we have to impose additional relations on $\BV^i_j$. The first one is that at each point $p=(\chi,\nu)\in B(\XX,\VV)(x)$ the relation:
\begin{equation}\label{warunek 0}
\pi_*(\BV^0_j(p))=V_j
\end{equation}
holds for every $j=1,\ldots,m$, where $\nu=(V_1,\ldots,V_m)$ is a basis of $\VV(x)$. Vector fields $\BV^i_j$ for $i\geq 1$ are defined by the relation:
\begin{equation}\label{warunek 1}
\BV^i_j=\ad^i_\BX\BV^0_j,
\end{equation}
where $\BX$ is the canonical vector field defined in Lemma \ref{lemat 4} and $j=1,\ldots,m$.  Conditions \eqref{warunek 0} and \eqref{warunek 1} define $\BV^i$ uniquely up to $\BG$ and $\BF$. In order to normalise $\BV^0$ in the vertical directions let us define functions $C_{pql}^{ir}$ by the equations:
$$
[\BV^0_p,\BV^i_q]=\sum_{l,r}C_{pql}^{ir}\BV^l_r\mod\BX,\BG,\BF.
$$
We introduce the following conditions:
\begin{equation}\label{warunek 2}
C_{pq1}^{1r}=0,
\end{equation}
for arbitrary $p,q,r=1,\ldots,m$,
\begin{equation}\label{warunek 3}
\sum_{p=1}^m C_{pq0}^{1p}=0,
\end{equation}
for arbitrary $q=1,\ldots,m$,
\begin{equation}\label{warunek 4}
\sum_{q=1}^m C_{pq2}^{2q}=0,
\end{equation}
for arbitrary $p=1,\ldots,m$,
\begin{equation}\label{warunek 4b}
\sum_{q=1}^m C_{pq3}^{3q}=0,
\end{equation}
for arbitrary $p=1,\ldots,m$.

Condition \eqref{warunek 2} will be responsible for the normalisation in the directions of $\BG$. Conditions \eqref{warunek 3} and \eqref{warunek 4} will be responsible for the normalisation in the directions of $\BF^0$ and $\BF^1$ in the case $m>1$. Conditions \eqref{warunek 3} and \eqref{warunek 4b} will be responsible for the normalisation in the directions of $\BF^0$ and $\BF^1$ in the case $k>2$ and $m=1$. Note that we use $\BV^3$ in \eqref{warunek 4b}. This is the point where $k>2$ is necessary for our constructions.
\begin{remark}
Assume that $m=1$. Then $\BV^i=(\BV^i_1)$ and by Jacobi identity we get
$$
[\BX,[\BV^0_1,\BV^1_1]] = [\BV^0_1,\BV^2_1]-[\BV^1_1,\BV^1_1] =[\BV^0_1,\BV^2_1].
$$
It follows that brackets $[\BV^0_1,\BV^1_1]$ and $[\BV^0_1,\BV^2_1]$ are related. This is the reason why we use $[\BV^0_1,\BV^3_1]$ instead of $[\BV^0_1,\BV^2_1]$ in the condition \eqref{warunek 4b}.
\end{remark}

\begin{theorem}\label{twierdzenie 2}
Let $(\XX,\VV)$ be a regular pair on a manifold $M$ of dimension $n=mk+m+1$, where $\rk\,\VV=m$. If $m>1$ and $k>1$ then there exists the unique frame on $B(\XX,\VV)$ satisfying conditions \eqref{warunek 0}-\eqref{warunek 4}. If $m=1$ and $k>2$ then there exists the unique frame on $B(\XX,\VV)$ satisfying conditions \eqref{warunek 0}-\eqref{warunek 3} and \eqref{warunek 4b}. Two pairs $(\XX,\VV)$ and $(\YY,\WW)$ are equivalent if and only if the corresponding frames on $B(\XX,\VV)$ and $B(\YY,\WW)$ are diffeomorphic.
The symmetry group of $(\XX,\VV)$ is at most $(n+m^2+2)$-dimensional and it is has maximal dimension if and only if $(\XX,\VV)$ is locally of equation type and the corresponding system of equations is trivial: $x^{(k+1)}=0$.  The following structural equations are satisfied:
\begin{eqnarray}\label{strukturalne 1}
&&[\BG^p_q,\BX]=0,\\
&&[\BF^0,\BX]=-\BX,\nonumber\\
&&[\BF^1,\BX]=-2\BF^0-k\sum_j\BG^j_j.\nonumber
\end{eqnarray}
Additionally, if $(\XX,\VV)$ is of equation type then:
\begin{eqnarray}\label{strukturalne 2}
&&[\BG^p_q,\BV^i_j]=\delta^p_j\BV^i_q,\\
&&[\BF^0,\BV^i_j]=-i\BV^i_j,\nonumber\\
&&[\BF^1,\BV^i_j]=i(i-1-k)\BV^{i-1}_j\nonumber.
\end{eqnarray}
\end{theorem}
\begin{proof}
The proof is divided into three parts. At the beginning we will construct canonical frame, then we will consider the most symmetric case and finally we will show that structural equations are satisfied.

{\bf Construction of the canonical frame.}
Let us fix a projective vector field $X$ and a normal frame $V$ of $\VV$. Then $V, \ad_XV,\ldots,\ad^k_XV, X$ is a frame on $M$. The group $\mathrm{GL}(m)$ acts on $V$ by multiplication. Assume that the condition \eqref{warunek 0} holds. In local coordinates on $B(\XX,\VV)$ we have:
$$
\BV^0_j=\sum_{p=1}^mG_j^pV_p+\langle\beta_j,\BG\rangle+ \gamma_{j0}\BF^0+\gamma_{j1}\BF^1,
$$
where $j=1,\ldots,m$, $\beta_j=(\beta^s_{jt})_{s,t=1,\ldots,m}$ and we use the abbreviation $\langle\beta_j,\BG\rangle =\tr(\beta_j\BG)= \sum_{s,t=1}^m\beta^s_{jt}\BG^t_s$. Our aim is to find functions $\beta_j,\gamma_{j0},\gamma_{j1}$, such that the conditions \eqref{warunek 2}-\eqref{warunek 4} (or \eqref{warunek 4b}) are satisfied. If we do this then the condition \eqref{warunek 1} will define the canonical frame on $B(\XX,\VV)$.

We use Lemma \ref{lemat 5} and compute that: 
\begin{eqnarray*}
\BV^1_j&=&[\BX,\BV^0_j]\\
&=&\sum_pG_j^p\frac{1}{F_0}(\ad_XV_p-kF_1V_p)\\
&+&\langle\BX(\beta_j),\BG\rangle +k\frac{1}{F_0}(\gamma_{j1}F_0-\gamma_{j0}F_1)\sum_p\BG_p^p +\gamma_{j0}\frac{1}{F_0}X \mod \BF,
\end{eqnarray*}
We use here a simple fact that $[\BG_s^t,\sum_j\BG_j^j]=0$. The next Lie brackets immediately give:
\begin{eqnarray*}
\BV^i_j&=&\sum_pG_j^p\frac{1}{(F_0)^i}\ad^i_XV_p\mod V,\ad_XV,\ldots,\ad_X^{i-1}V,\BX,\BG,\BF,
\end{eqnarray*}
for $i=2,\ldots,k$. More precisely we have:
\begin{eqnarray*}
\BV^i_j&=&\sum_pG_j^p\frac{1}{(F_0)^i}(\ad^i_XV_p +c^i_1F_1\ad^{i-1}_XV +\cdots\\&+&c^i_{i-1}(F_1)^{i-1}\ad_XV +c^i_i(F_1)^iV) \mod\BX,\BG,\BF,
\end{eqnarray*}
where $c^i_j$ are certain rational numbers (their exact values are not important for us). By direct calculations we find that:
\begin{eqnarray*}
[\BV^0_p, \BV^1_q]&=&\sum_{s,t} G_p^sG_q^t\frac{1}{F_0}([V_s,\ad_XV_t]-kF_1[V_s,V_t])\\
&+&\sum_{s,r}\beta_{pq}^sG_s^r\frac{1}{F_0}(\ad_XV_r-kF_1V_r)
-\sum_{s,r}\BX(\beta_{qp}^s)G_s^rV_r\\
&-&k\sum_r\frac{1}{F_0}(\gamma_{p1}G_q^rF_0-\gamma_{p0}G^r_qF_1 +\gamma_{q1}G_p^rF_0-\gamma_{q0}G^r_pF_1)V_r\\
&-&\sum_r(\gamma_{p0}G_q^r+\gamma_{q0}G_p^r)\frac{1}{F_0}\ad_XV_r \mod\BX,\BG,\BF,
\end{eqnarray*}
\begin{eqnarray*}
[\BV^0_p, \BV^2_q]&=&\sum_{s,t}G_p^sG_q^t\frac{1}{(F_0)^2}\left( \sum_{l=0}^2c^2_l(F_1)^l[V_s,\ad^{2-l}_XV_t]\right)\\
&+&\sum_{s,r}\beta_{pq}^sG_s^r\frac{1}{(F_0)^2}\ad^2_XV_r\\
&-&2\sum_r\gamma_{p0}G_q^r\frac{1}{(F_0)^2}\ad_X^2V_r \mod\BV^0,\BV^1,\BX,\BG,\BF,
\end{eqnarray*}
and:
\begin{eqnarray*}
[\BV^0_p, \BV^3_q]&=&\sum_{s,t}G_p^sG_q^t\frac{1}{(F_0)^3}\left( \sum_{l=0}^3c^3_l(F_1)^l[V_s,\ad^{3-l}_XV_t]\right)\\
&+&\sum_{s,r}\beta_{pq}^sG_s^r\frac{1}{(F_0)^3}\ad^3_XV_r\\
&-&3\sum_r\gamma_{p0}G_q^r\frac{1}{(F_0)^3}\ad_X^3V_r \mod\BV^0,\BV^1,\BV^2,\BX,\BG,\BF.
\end{eqnarray*}
Equivalently we can write:
\begin{eqnarray*}
[\BV^0_p, \BV^1_q]&=&\sum_{s,t} G_p^sG_q^t\frac{1}{F_0}([V_s,\ad_XV_t]-kF_1[V_s,V_t])\\
&-&\sum_r\BX(\beta_{qp}^r)\BV^0_r -k\gamma_{p1}\BV^0_q-k\gamma_{q1}\BV^0_p\\
&+&\sum_r\beta_{pq}^r\BV^1_r -\gamma_{p0}\BV^1_q-\gamma_{q0}\BV^1_p\mod\BX,\BG,\BF,
\end{eqnarray*}
\begin{eqnarray*}
[\BV^0_p, \BV^2_q]&=&\sum_{s,t}G_p^sG_q^t\frac{1}{(F_0)^2}\left( \sum_{l=0}^2c^2_l(F_1)^l[V_s,\ad^{2-l}_XV_t]\right)\\
&+&\sum_r\beta_{pq}^r\BV^2_r -2\gamma_{p0}\BV^2_q\mod\BV^0,\BV^1,\BX,\BG,\BF,
\end{eqnarray*}
and:
\begin{eqnarray*}
[\BV^0_p, \BV^3_q]&=&\sum_{s,t}G_p^sG_q^t\frac{1}{(F_0)^3}\left( \sum_{l=0}^3c^3_l(F_1)^l[V_s,\ad^{3-l}_XV_t]\right)\\
&+&\sum_r\beta_{pq}^r\BV^3_r -3\gamma_{p0}\BV^3_q\mod\BV^0,\BV^1,\BV^2,\BX,\BG,\BF.
\end{eqnarray*}

Since vector fields $\ad^iV_j$ together with $X$ span the whole tangent bundle $TM$, we can express the Lie bracket $[V_s,\ad^l_XV_t]$ as a linear combination of $\ad^i_XV_j$ and $X$. If we lift them to $B(\XX,\VV)$, we can also write $[V_s,\ad^l_XV_t]=\sum_{i,j}c_{sti}^{lj}\BV^i_j\mod \BX,\BG,\BF$, for some functions $c_{sti}^{lj}$ on $B(\XX,\VV)$. Then we are able to rewrite the conditions \eqref{warunek 2}-\eqref{warunek 4} in terms of unknown functions  $\beta_{js}^t$, $\gamma_{j0}$ and $\gamma_{j1}$. We get the following system of equations:
\begin{eqnarray}
\beta_{pq}^r-\delta^r_p\gamma_{p0}-\delta^r_q\gamma_{q0} &=&C_1(p,q,r),\label{uklad 1}\\
-\sum_p\BX(\beta_{qp}^p)-k(m+1)\gamma_{q1} &=&C_2(q),\label{uklad 2}\\
\sum_q\beta_{pq}^q-2m\gamma_{p0} &=&C_3(p),\label{uklad 3}
\end{eqnarray}
where $C_1$, $C_2$, $C_3$ are certain functions on $B(\XX,\VV)$. If $m=1$ and $k>2$ the last equation is replaced by:
\begin{equation}\label{uklad 4}
\sum_q\beta_{pq}^q-3m\gamma_{p0}=C_4(p).
\end{equation}
Note that the equations are homogeneous of order one in $G$. They can be solved and the solution is unique. We proceed as follows. In the first step we find $\gamma_{p0}$ using \eqref{uklad 1} and \eqref{uklad 3} or \eqref{uklad 4} (depending whether $m>1$ or $m=1$). Namely, we substitute $r=q$ in \eqref{uklad 1}, then we take the sum $\sum_q$ and subtract the result from \eqref{uklad 3} or \eqref{uklad 4} in order to eliminate $\sum_q\beta^q_{pq}$. We get a system of equations for $\gamma_{p0}$, which has the unique solution. In this way we get $\gamma_{p0}$, which we substitute to \eqref{uklad 1} and we find $\beta_{pq}^r$. Finally, from \eqref{uklad 2} we get $\gamma_{q1}$.

From the form of the equations we can deduce that $\beta_{ps}^t$, $\gamma_{r0}$ and $\gamma_{r1}$ are homogeneous of order one in $G^s_t$.

Now we can also see why we need the condition \eqref{warunek 4b} which involves $\BV^3$ in the case $m=1$. In this case the left hand sides of \eqref{uklad 1} and \eqref{uklad 3} coincide and thus we can not use \eqref{warunek 4} (compare Remark which precedes the theorem).

{\bf Uniqueness of the model with maximal symmetry group.}
If the dimension of the symmetry group of $(\XX,\VV)$ is maximal, equal $\dim B(\XX,\VV)$, then all structural functions of the canonical frame have to be constant. At the beginning we use Lemma \ref{lemat 5} and directly compute:
$$
[\BF^0,\BX]=-\BX,\qquad [\BF^1,\BX]=-2\BF^0-k\sum_j\BG^j_j.
$$
We also find that:
$$
[\BG^p_q,\BX]=0.
$$
In particular we proved \eqref{strukturalne 1}.

Note that if a certain non-vanishing structural function is constant then it has to be homogeneous of order 0 in any coordinate. In our case the functions $\beta_j, \gamma_{j0}, \gamma_{j1}$ are homogeneous of order one in $G=(G^i_j)$. This implies that also vector fields $\BV^i_j$ are homogeneous of order one in $G$. On the other hand $\BX,\BF,\BG$ are homogeneous of order 0 in $G$. Therefore, in the most symmetric case, we get that
$[\BG^p_q,\BV^i_j]=0\mod\BV$. Precisely we have:
$$
[\BG^p_q,\BV^i_j]=\delta^p_j\BV^i_q,
$$
as can be easily seen from the formula for $\BV^i_j$. The same argument applies to $[\BF^0,\BV^0]$ and $[\BF^1,\BV^0]$. We have:
$$
[\BF^0,\BV^0]=[\BF^1,\BV^0]=0.
$$
If $[\BF^0,\BV^0]=[\BF^1,\BV^0]=0$ then also structural functions of brackets $[\BF^0,\BV^i]$ and $[\BF^1,\BV^i]$ are constant. They can be computed using Jacobi identity and equation \eqref{warunek 1}. By homogeneity argument we also get:
$$
[\BV^i,\BV^l]=0.
$$
By definition:
$$
[\BX, \BV^i]=\BV^{i+1},
$$
for $i=0,\ldots,k-1$ and there remains the last bracket $[\BX, \BV^k]$ to find. We have:
$$
[\BX,\BV^k]= w_0\BV^0+\cdots+w_{k}\BV^{k},
$$
for some $w_i$. We claim that all $w_i$ are polynomial of order $k+1-i$ in $F_0$. Indeed in local coordinates on $B(\XX,\VV)$ the vector fields $\BV^i$ are homogeneous of order $-i$ in $F_0$. Additionally we see that $[\BX,\BV^k]$ is of order $-(k+1)$. Thus $w_i$ has to be of order $k+1-i$. As a conclusion we get that all $w_i$ vanish provided that they are constant.
In this way we have proved uniqueness of the most symmetric model of pairs $(\XX,\VV)$. Therefore it has to be locally equivalent to the pair corresponding to the trivial system $x^{(k+1)}=0$.

{\bf Structural equations.}
We have already proved that equations \eqref{strukturalne 1} are satisfied. We will show now \eqref{strukturalne 2}. We easily compute in coordinates that:
\begin{eqnarray}\label{strukturalne 2a}
&&[\BG^p_q,\BV^0_j]=\delta^p_j\BV^0_q + \sum_{s,t}P^{pt}_{qjs}\BG^s_t +Q^p_{qj}\BF^0+R^p_{qj}\BF^1,\\
&&[\BF^0,\BV^0_j]=\sum_{s,t}P^{0t}_{js}\BG^s_t +Q^0_j\BF^0+R^0_j\BF^1,\nonumber\\
&&[\BF^1,\BV^0_j]=\sum_{s,t}P^{1t}_{js}\BG^s_t +Q^1_j\BF^0+R^1_j\BF^1,\nonumber
\end{eqnarray}
for some coefficients $P$, $Q$ and $R$. The Jacobi identity applied to $[\BX,[\BF^0,\BV^i_j]]$ and $[\BX,[\BF^1,\BV^i_j]]$ gives the formulae for $[\BF^0,\BV^{i+1}_j]$ and $[\BF^1,\BV^{i+1}_j]$. Indeed, by induction we prove (we use \eqref{warunek 1} and \eqref{strukturalne 1}):
\begin{eqnarray}\label{strukturalne 2b}
&&[\BG^p_q,\BV^i_j]=\delta^p_j\BV^i_q\mod \BX,\BG,\BF,\\
&&[\BF^0,\BV^i_j]=-i\BV^i_j\mod \BX,\BG,\BF,\nonumber\\
&&[\BF^1,\BV^i_j]=i(i-1-k)\BV^{i-1}_j\mod \BX,\BG,\BF.\nonumber
\end{eqnarray}
Moreover, if $P=Q=R=0$ then \eqref{strukturalne 2} hold for every $i$. Therefore in order to finish the prove we only have to show that the vertical components vanish in \eqref{strukturalne 2a}.

Let us notice that if $(\XX,\VV)$ is of equation type then:
\begin{equation}\label{relacja 0}
[\BV^0_p,\BV^i_q]=0\mod \BX,\BG,\BF,\BV^0,\ldots,\BV^i,
\end{equation}
i.e. there is no component next to $\BV^l$ for $l>i$. This is a consequence of condition (G4), since if (G4) holds then projection of $\BV^0_p$ to the base manifold is contained in the Cauchy characteristic of any distribution $\VV^i$. Moreover, if $(\XX,\VV)$ is of equation type then condition \eqref{warunek 2} implies:
\begin{equation}\label{relacja 1}
[\BV^0_p,\BV^0_q]=0\mod \BX,\BG,\BF.
\end{equation}
If not then Jacobi identity applied to $[\BX,[\BV^0_p,\BV^0_q]]$ would contradicts \eqref{warunek 2}.

Now, we can consider all expressions of the following type:
$$
[\BG^p_q,[\BV^0_u,\BV^i_w]],\qquad [\BF^0,[\BV^0_u,\BV^i_w]],\qquad [\BF^1,[\BV^0_u,\BV^i_w]],
$$
where $i=1,2$ when $m>1$ or $i=1,3$ when $m=1$. If we apply, once again, Jacobi identity and take into account conditions \eqref{warunek 2}-\eqref{warunek 4} (or \eqref{warunek 4b}) we obtain a system of linear equations for functions $P$, $Q$ and $R$. If \eqref{relacja 0} and \eqref{relacja 1} hold then the unique solution has the form $P=Q=R=0$.
\end{proof}
\begin{remark}
Our main Theorem \ref{twierdzenie 0} is a simplified version of Theorem \ref{twierdzenie 2}.
\end{remark}

\section{Cartan Connection}
Assume that $L$ is a Lie group and $L_0$ is a closed Lie subgroup of $L$. Denote by $\fl$ and $\fl_0$ the corresponding Lie algebras. Let $M$ be a manifold of dimension $n=\dim L/L_0$ and $P$ be a principal bundle over $M$ with group $L_0$. We say that one-form $\hat\omega$ on $P$ with values in $\fl$ is \emph{Cartan connection of type $(L,L_0)$} if
\begin{enumerate}
\item $\hat\omega(l^*)=l$ for every $l\in\fl_0$, where $l^*$ is a fundamental vector field on $P$ defined by $l$,
\item $(R_H)^*\hat\omega=Ad_{H^{-1}}\hat\omega$ for every $H\in L_0$, where $Ad_{H^{-1}}$ is the adjoint action of $H^{-1}$ on $\fl$,
\item $\hat\omega\colon T_xP\to\fl$ is an isomorphism for every $x\in P$.
\end{enumerate}

In our situation $B(\XX,\VV)$ is principal $\mathrm{GL}\oplus \mathrm{PGL}_0(m)$-bundle. Therefore $L_0=\mathrm{GL}\oplus \mathrm{PGL}_0(m)$ and we want to find $L$. Let us denote
$$
\BH=2\BF^0+k\sum_j\BG^j_j,\qquad \BY=\BF^1,
$$
and
$$
\BW^i_j=\frac{1}{i!}\BV^i_j
$$
\begin{proposition}\label{stwierdzenie 1}
Let $(\XX,\VV)$ be a regular pair. The following structural equations are satisfied:
\begin{equation}\label{strukturalne 3}
[\BX,\BY]=\BH,\qquad [\BH,\BX]=-2\BX,\qquad [\BH,\BY]=2\BY,
\end{equation}
\begin{equation}
[\BG,\BX]=0,\qquad [\BG,\BY]=0,\qquad [\BG,\BH]=0.\nonumber
\end{equation}
Additionally, if $(\XX,\VV)$ is of equation type then:
\begin{equation}
[\BX,\BW^i_j]=(i+1)\BW^{i+1}_j,\quad [\BY,\BW^i_j]=-(k-i+1)\BW^{i-1}_j,\quad
[\BH,\BW^i_j]=-2i\BW^i_j,\nonumber
\end{equation}
\begin{equation}\label{strukturalne 4}
[\BG^p_q,\BW^i_j]=\delta^p_j\BW^i_q.
\end{equation}
\end{proposition}
\begin{proof}
Follows directly from the structural equations \eqref{strukturalne 0}, \eqref{strukturalne 1}, \eqref{strukturalne 2} and definitions of $\BH$, $\BY$ and $\BW^i_j$.
\end{proof}

Let 
$$
(\alpha,\beta,\gamma,\theta_0^j,\ldots,\theta_k^j,\omega^t_s\ |\ j,s,t=1,\ldots,m),
$$
be the coframe on $B(\XX,\VV)$ dual to the frame
$$
(\BH,\BY,\BX,\BV^0_j,\ldots,\BV^k_j,\BG^s_t\ |\ j,s,t=1,\ldots,m).
$$
Let us briefly write
$$
\theta_i=(\theta_i^1,\ldots,\theta_i^m)^t, \qquad \theta=(\theta_0,\ldots\theta_k), \qquad \omega=(\omega^t_s)_{s,t=1,\ldots,m},
$$
and introduce
$$
\tilde\omega=
\left(
\begin{array}{cc}
0&0\\
\theta&\omega
\end{array}
\right).
$$
$\tilde\omega$ is a 1-form with values in the Lie algebra
$$
\mathfrak{a}=\left\{\left(
\begin{array}{cc}
0&0\\
v&g
\end{array}
\right)\in \mathfrak{gl}(m+k+1)\ |\ v\in\R^{m\times(k+1)}, g\in\mathfrak{gl}(m)\right\}
$$
of the following Lie group of matrices 
$$
A=\left\{\left(
\begin{array}{cc}
\Id_{k+1}&0\\
V&G
\end{array}
\right)\in \mathrm{GL}(m+k+1)\ |\ V\in\R^{m\times(k+1)}, G\in\mathrm{GL}(m)\right\}.
$$
One can treat $A$ as an "extended" affine group. $\mathfrak{gl}(m)$ is a Lie sub-algebra of $\mathfrak{a}$.

Let us also introduce the form
$$
\eta=\left(\begin{array}{cc}
\alpha&\beta\\
\gamma&-\alpha\\
\end{array}\right).
$$
with values in Lie algebra $\mathfrak{sl}(2)$. Note that $\pgl(2)_0$ can be embedded into $\mathfrak{sl}(2)$ in the following way:
$$
\pgl(2)_0=\left\{\left(
\begin{array}{cc}
\frac{1}{2}a&b\\
0&-\frac{1}{2}a\\
\end{array}
\right)\ |\ a,b\in\R\right\}.
$$
\begin{theorem}\label{twierdzenie 3}
Assume that a pair $(\XX,\VV)$ is of equation type. The form $\hat\omega=(\eta,\tilde\omega)$ is Cartan connection of type $(L, \mathrm{PGL}(2)_0\oplus\mathrm{GL}(m))$, where $L$ is a semi-direct product of $\mathrm{SL}(2)$ and $A$, where $\mathrm{SL}(2)$ acts irreducibiy on each copy of $\R^{k+1}$.
\end{theorem}
\begin{proof}
Follows directly from Proposition \ref{stwierdzenie 1}.
\end{proof}

\end{document}